\documentclass{gtart}
\usepackage{amssymb,amsmath,euscript,rlepsf,latexsym}
\usepackage[all]{xy}

\input gtspec
\volumenumber{3}\papernumber{10}\volumeyear{1999}
\pagenumbers{235}{252}\published{2 September 1999}
\proposed{Robion Kirby}
\seconded{Cameron Gordon, Joan Birman}
\received{21 May 1999}\revised{27 August 1999}
\accepted{1 September 1999}

\newtheorem{thm}{Theorem}
\newtheorem{lem}{Lemma}
\newtheorem{cor}[lem]{Corollary}

\newtheorem*{conj}{Conjecture}

\theoremstyle{definition}

\theoremstyle{remark}

\newtheorem{rem}{Remark}

\newcommand{\Ker}{\operatorname{Ker}}

\newcommand{\ind}{\operatorname{index}}

\newcommand{\Z}{\mathbb{Z}}
\newcommand{\N}{\mathbb{N}}
\newcommand{\R}{\mathbb{R}}

\newcommand{\C}{\mathbb{C}}
\newcommand{\e}{\epsilon}
\newcommand{\hra}{\hookrightarrow}
\newcommand{\imra}{\looparrowright}

\newcommand{\ra}{\longrightarrow}

\newcommand{\lra}{\leftrightarrow}

\newcommand{\svh}{\Sigma_\nu H}
\newcommand{\svf}{\Sigma_\nu F}

\newcommand{\svhp}{\Sigma_\nu^+ H}
\newcommand{\svfp}{\Sigma_\nu^+ F}

\newcommand{\om}{\widetilde{\Omega}}
\newcommand{\x}{\times}
\newcommand{\dd}{\partial}

\begin{document}

\title{All two dimensional links are null homotopic}
\authors{Arthur Bartels\\Peter Teichner}

\address{University of California in San Diego\\La Jolla, CA,
92093-0112, USA}

\email{abartels@math.ucsd.edu\\teichner@math.ucsd.edu}

\asciiabstract{We show that any number of disjointly embedded 
2-spheres in
4-space can be pulled apart by a link homotopy,
ie, by a motion in which the 2-spheres stay disjoint but are
allowed to self-intersect.}

\begin{abstract}
We show that any number of disjointly embedded $2$--spheres in
$4$--space can be pulled apart by a {\em link homotopy},
ie, by a motion in which the $2$--spheres stay disjoint but are
allowed to self-intersect.
\end{abstract}



\maketitlepage

\section{Introduction}

In order to separate $3$--dimensional linking and knotting
phenomena, John Milnor
introduced in 1954 the notion of a {\em link homotopy} \cite{M}. It is a
one-parameter family of maps
$$
S^1 \amalg \dots \amalg S^1 \to \R^3
$$
 during which one allows self-intersections but does not allow
different components to cross. For example, any knot is link
homotopically trivial but
the Hopf link is not. Another example of a homotopically essential
link in $\R^3$ is given by the
Borromean rings which are detected by the generalized linking number
$\mu(1,2,3)$. More generally, Milnor showed that a link in $\R^3$ is
homotopically trivial if and only if all its $\mu$--invariants (with
non-repeating indices) vanish. Another way to formulate this result is to
consider a certain quotient of the fundamental group of the link complement,
now known as the {\em Milnor group}, which is an invariant of link homotopy
(see Section~\ref{sec:groups}). Then a link is homotopically trivial
if and only if its Milnor group is isomorphic to the Milnor group
$MF_n$ of the unlink (with $n$ components).
In our paper this particular group $MF_n$, the {\em free} Milnor
group, is used as the key ingredient to prove the following result:

\begin{thm} \label{main}
Every smooth link $L\co  S^2 \amalg \dots \amalg S^2 \hra \R^4$ is link
homotopic to the unlink.
\end{thm}
The beginning of 4--dimensional link homotopy was the paper
\cite{FR} by Roger Fenn and Dale Rolfsen in 1985 who construct two disjointly
immersed $2$--spheres in $\R^4$ which are not link homotopically trivial.
William Massey and Rolfsen had just introduced higher dimensional link
homotopy and observed that their generalized linking number for two
$2$--spheres in $\R^4$ vanishes on embedded links. They ask in
\cite{MR} whether
Theorem~\ref{main} is true for two component links, hence this question is
sometimes referred to as the {\em Massey--Rolfsen problem}. A proof has been
attempted several times but as the referee points out, this paper gives the
first correct solution, as well as a generalization to arbitrary many
components.

In the course of our proof we have to introduce many
self-intersections into the components of $L$ but surprisingly we
can keep different components disjoint.
The argument has two
completely independent steps: one is to construct a link concordance to the
unlink and the other is to improve the link concordance to a link homotopy.

Both of these steps generalize to links of $n$--spheres in
$\R^{n+2}$ for all $n>1$. The first step is Bartels' PhD thesis
\cite{Ba}, the second Teichner's habilitation \cite{habil}. These papers
are long and yet unpublished whereas both steps are discussed in
full detail in  this short note in dimension~4. We would still like to announce
the general result, since it seems to come as a surprise how far the
Massey--Rolfsen problem can be pushed. For readers interested in high
dimensional link homotopy we should mention that the theory was developed
by Fenn, Habegger, Hilton, Kaiser, Kirk, Koschorke, Massey, Nezhinsky, Rolfsen
and others.

\begin{thm}
For $n>1$, every smooth link $L\co  S^n \amalg \dots \amalg S^n \hra \R^{n+2}$ is
link homotopic to the unlink.
\end{thm}

The result for two $2$--spheres in $S^4$ is also proven by very different
methods in \cite{T}. This paper, which is currently being rewritten, actually
gives a complete calculation of the group $LM_{2,2}^4$ of link homotopy
classes of link maps $S^2\amalg S^2
\to\R^4$. In particular,  it implies the two component case of the following
conjecture.
\begin{conj}
Theorem~\ref{main} still holds if one component of the link $L$ is
not embedded (but mapped into $\R^4$ disjointly from the other
components).
\end{conj}
The conjecture is supported by the fact that it holds for one-- and
two--compon\-ent links and that all known invariants
vanish on links with only one non-embedded component. It is not difficult
to construct three disjointly immersed $2$--spheres in $\R^4$, one of
them embedded, which are not homotopically trivial.

Here is a brief outline of the proof for Theorem~\ref{main}.
The last section (Section~\ref{sec:habil}) contains the proof that ``link
concordance implies link homotopy''   for link maps
in dimension~4. The main idea (which works for arbitrary dimensions as long
as the {\em codimension} is $\geq 2$) is to develop a theory of ambient {\em
singular} handles. This leads to a generalization to immersions of Colin
Rourke's ambient handle proof
\cite{R} of Hudson's theorem that ``concordance implies isotopy'' in
codimension $\geq 3$. The proof given here simplifies in contrast to higher
dimensions out of several reasons:
\begin{itemize}
\item there are no triple points,
\item the singular handles are given by the well known 4--dimensional Whitney
move (and it's reverse),
\item there is no need to use ambient handle slides (as in \cite{R}) or ambient
Cerf theory (as in
\cite{habil}) because the relevant product structures are guaranteed by simple
arguments involving only $0$--handles.
\end{itemize}

In the other sections we show that the link
$L$ is link concordant to the unlink, ie, that it bounds disjoint
immersions of
$3$--balls into
$D^5$. One observes that by possibly  introducing more singularities (finger
moves) into the $3$--balls one may assume that the fundamental group of the
complement of these $3$--balls is the free Milnor group
$MF_n$. Therefore, we  first construct a certain $5$--manifold with
fundamental group
$MF_n$ which plays the role of the complement of the
$3$--balls in $D^5$ and we then fill it in with standard thickenings of immersed
$3$--balls to get back $D^5$. The first step is to ask whether $0$--surgery on
$L$ bounds a
$5$--manifold over the group $MF_n$. The answer to this question is ``yes''
using
the following two known algebraic facts \cite{M}, \cite{BK}
(which will be explained in the next section):
\begin{itemize}
\item $MF_n$ is a finitely presented nilpotent group.
\item The tower of nilpotent quotients of the free group $F_n$ is
homologically pro-trivial.
\end{itemize} By Alexander duality and Stallings' theorem \cite{S}, the
nilpotent quotients of $\pi_1(S^4  -  L)$ are the same as of $F_n$ which will
prove the existence of the searched for $5$--manifold. The argument finishes
with some surgeries which correct the second homology such that the filling
actually leads to a $5$--manifold which is recognized to be the $5$--ball by an
application of the
$6$--dimensional h--cobordism theorem.

\rk{Acknowledgments}  We wish to thank Mike Freedman and Colin Rourke for
interesting comments.

The second author is supported by NSF Grant DMS 97-03996.

\section{Milnor groups} \label{sec:groups}
We first collect the necessary group theoretic
facts:  The {\em lower central series} of a group
$G$  is defined by $G_1:= G, G_{k+1}:= [G,G_k]$
for $k\ge 1$.

All nilpotent quotients of a group factor through some
$G/G_k$, however these are not the only interesting
nilpotent quotients.  If $G$ is provided with a set of
normal generators, $G=\ll\! x_1, \cdots, x_n\!\gg$, then
one defines the {\it Milnor group}
$$
MG:=G/\ll\![ x_i, x^y_i]\!\gg,\ \ i=1, \dots, n \ \text{and}\
\  y \in G.
$$
These quotients were introduced in Milnor's
thesis \cite{M}.
The following lemma about Milnor
groups uses the fact that we have chosen the
elements $x_{i}$ to be normal generators for $G$.

\begin{lem}
\label{nilpotent}$MG$ has nilpotency class
$\le n$,
ie $MG_{n+1}=\{1\}$.
\end{lem}

\begin{proof}  We will use an
induction on the number $n$  of normal generators.  If $n=1$
the relations imply that all $x^y_{i}$  commute in
$MG$.  Since by assumption these elements generate $G$ it
follows that all commutators vanish in $MG$.

Now   assume the statement
holds for groups with $n-1$ normal generators and let $G$ be
normally generated by $x_{i},\ 1 \le i \le n$.   Define $A_i
\trianglelefteq MG$ to be
the normal closure of the element $x_{i}$.  Since all the
conjugates of this element
commute, $A_i$ is abelian.  Moreover, the intersection
$A$ of all $A_i$ lies obviously in the center of $MG$.
Now consider a commutator $[x,y]$ with $x\in MG, y\in
MG_n$.  Since all quotients $MG/A_i$ are
Milnor groups with
$n-1$ generators it follows by induction that $y\in A$
and thus $y$ is central ie $[x,y]=1$ in $MG$.
This shows that $MG_{n+1} = \{1\}$.
\end{proof}

\begin{cor}
\label{finitely_presented}
$MG$ is generated by
$x_{1},\dots, x_n$ and is also finitely presented.
\end{cor}

\begin{proof}
The statement
for the generators follows from the standard {\em
rewriting process} in nilpotent groups:  If a nilpotent
group $N$ is normally generated by $x_i$ then it is also
generated by these elements.  One uses an induction on
the nilpotency class of $N$ based on the fact that
$x\equiv y$ mod $N_k$ implies $a^x\equiv a^y$ mod
$N_{k+1}$ for all $a, x, y \in N$.  Moreover, the
fact that $N_k$ is generated by $k$--fold
commutators $[x_{i_1}, \dots, x_{i_k}]$ if the
$x_i$ generate $N$ shows that $N_k$  is finitely
generated if $N$ is.  An
induction on the nilpotency class together with
the fact that a (central) extension of finitely
presented groups is finitely presented implies that
a finitely generated nilpotent group is also
finitely presented.
\end{proof}

We need one more result which is basically due to Bousfield and Kan
\cite{BK}, and also uses some work of Ferry \cite{F} and Cochran
\cite{C}. Let $\Omega$ be any generalized homology theory. For a
group $G$ we abbreviate
the reduced theory by
\[
\om_k(G):=\Ker (\Omega_k(K(G,1)) \to \Omega_k(*) )
\]
where $\pi_1K(G,1)=G$ and $\pi_iK(G,1)=0$ for $i>1$.
 \begin{thm} \label{eventual}
Let $F$ be the free group on $n$ generators.
Given $r\in\N$ and $k>1$ such that $\Omega_{k-1}(*)=0$, there
exists an integer $d=d(k,r,n)$
such that the map
$$
\om_k(F/F_{r+d}) \ra \om_k(F/F_r)
$$
is trivial.
\end{thm}
\begin{proof}
For $\Omega=$ ordinary homology, the theorem comes from \cite{BK}. To
get the general statement, one uses an {\em eventual Hurewicz
Theorem} as in \cite{C}.
It is therefore necessary to reduce to the case of simply-connected
spaces. This can be easily done in our context by picking maps
$f_i\co S^1 \to K(F/F_r,1)$, which represent the generators $x_i$ of
$F$, and
attach 2--cells to get simply-connected complexes $X_r$. Then
$$
H_k(X_r)\cong H_k(F/F_r)
$$
for all $k>1$, so the again by \cite{BK} the maps $H_k(X_{r+d})\to
H_k(X_r)$ are eventually zero.
By the eventual Hurewicz
Theorem for simply-connected spaces it follows that the maps
$X_{r+d}\to  X_r$ are eventually null homotopic and thus the induced maps
$\om_k(X_{r+d})\to \om_k(X_r)$ are eventually zero.
Moreover, there is an exact sequence
$$
\om_{k+1}(\vee^n S^2) \ra \om_k(F/F_r) \overset{i_*}{\ra} \om_k(X_r)
$$
By assumption and excision we have
$$
\om_{k+1}(\vee^n S^2) \cong \bigoplus^n \Omega_{k-1}(*) =0
$$
which proves that $i_*$ is a monomorphism. The same exact argument
works for $r$ replaced by $r+d$ which implies our claim.
\end{proof}

\section{Singular handles}

\begin{figure}[ht]
\vspace{.01in}
\hspace*{\fill}
\epsfysize=1in
\epsfbox{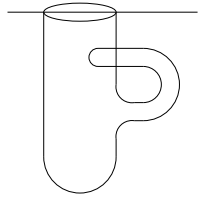}
\hspace*{\fill}
\vspace{.01in}
\caption{\label{singhand} The immersion $\Delta_1$}
\end{figure}

Let $\Delta_0 \co  D^3 \hra D^5$ be a standardly embedded slice
disk for an
unknot $S^2 \subset S^4$. Let $\Delta_\nu \co  D^3 \imra D^5$ be
obtained from $\Delta_0$ by
performing $\nu$ finger moves on $\Delta_0$ along arcs
$\gamma_i$ in the interior of
$D^5$ connecting pairs of points on $\Delta_0$. So $\Delta_\nu$
is an immersed slice disk for
the unknot, see figure~\ref{singhand}.
The self intersections of $\Delta_\nu$ consist of a
disjoint union of
$\nu$ circles.
Now extend $\Delta_\nu$ to a thickening
$\widetilde{\Delta_\nu} \co  D^3 \times D^2 \to D^5$
and define
\[
\begin{array}{lll}
\svh  & :=  & \widetilde{\Delta_\nu} ( D^3 \times D^2 )
\\
\svf  & :=  & \mbox{closure of } (\partial \svh - S^4)
\\
\svfp & :=  & \svf \cup_{S^2 \times S^1} D^3 \times S^1
\\
\svhp & :=  & \svh \cup_{S^2 \times D^2} D^3 \times D^2.
\\
\end{array}
\]
Note that
$\partial \svh = \svf \cup_{S^2 \times S^1} S^2 \times D^2 $
and
$\svfp = \partial \svhp$.
$\svh$ can also be constructed (abstractly) as a self-plumbing
of $D^3 \times D^2$ along the self intersection circles.
Similarly,
$\svhp$ is a self-plumbing of $S^3 \times D^2$.

If $P=\Delta_\nu (x)$, then
$m_P := \widetilde{\Delta_\nu} (x \times S^1)$ is called the
{\em meridian} to $\Delta_\nu$ at $P$.
Let $P = \Delta_\nu(x) = \Delta_\nu(y)$ be a self intersection
point of $\Delta_\nu$.
(So at $P$ there are actually two different meridians.)
Let $B \co  \R \times \C \times \C \hra D^5$ be a
parameterization of an open neighborhood
of $P$ such that
\[
\begin{array}{lllll}
B(0,0,0)              & = &  P
\\
B \cap \Delta_\nu     & = &  \R \times \C \times 0
 & \cup & \R \times 0   \times \C
\\
B \cap \svh           & = &  \R \times \C \times D^2
 & \cup & \R \times D^2 \times \C \mbox{ .}
\\
\end{array}
\]
Then $T_P := B(0 \times S^1 \times S^1) \subset \svf$ is called
the
{\em Clifford torus} at $P$. Note that $T_P$ bounds a solid torus,
say $B(0 \times S^1 \times D^2)$, in
$\svh$. Choose an arc in $D^3$ connecting $x$ to $y$. The image
under $\Delta_\nu$ of this arc gives a closed curve in $\svh$ and
a choice of an arc connecting this curve to a base point leads to an
element in $\pi_1 \svh$.
Applying this procedure to one point of each self-intersection
component gives
$\nu$ elements $g_1, \dots ,g_\nu \in \pi_1 \svh$, and it is
not hard to see that $\pi_1 \svh$ is
the free group on $\nu$ generators $g_i$. The $g_i$ have preimages in
$\pi_1 \svf$ which we will again denote by $g_i$. Let
$m \in \pi_1 \svf$ be the element induced by a fixed meridian
(to the unknot $S^2 \subset S^4$ say). The existence of the
Clifford tori in
$\svf$ implies the equations
\[
[m,m^{g_i}] = 1 \in \pi_1 \svf
\]
 for $i=1, \dots ,\nu$.

\begin{lem}
\label{onto}
$H_2 \svf$ maps onto $H_2 \svh$ and the Clifford tori
represent elements in the kernel of this map.
\end{lem}

\begin{proof}
The Mayer--Vietoris sequence
\[
0=H_3 D^5 \to H_2 \svf \to
H_2 \svh \oplus H_2 (D^5- \svh) \to
H_2 D^5=0.
\]
implies the first statement. The Clifford tori represent
trivial elements in $H_2 \svh$, since they bound solid tori there.
\end{proof}

Note that $\pi_1 \svf \cong \pi_1 \svfp$ and
$\pi_1 \svh \cong \pi_1 \svhp$. The meridian is an
element of the kernel of
$\pi_1 \svfp \to \pi_1 \svhp$. We will now construct a second
manifold $W$ bounding $\svfp$ such that $m$ and $g_1, \dots g_\nu$
give nontrivial elements in $\pi_1 W$.
Let $K_0 \co  S^3 \hra S^5$ be an unknot. Let
$E:= S^5 \#^\nu S^1 \times S^4$ where the connected sum operation
is done
away from $K_0$. So we still have
$K_0 \co  S^3 \hra E$. Now perform
finger moves on $K_0$ along $\nu$ arcs $\beta_1, \dots ,\beta_\nu$,
each following one of the handles just added, to obtain
$K_\nu \co  S^3 \imra E$. Thicken this to
$\widetilde{K_\nu}\co  S^3 \times D^2 \imra E$.
Note that $\widetilde{K_\nu} (S^3 \times D^2)$ is just $\svhp$. Now
$W:=E - \mbox{ Interior of } \svhp$ is the desired manifold.
$W$ is a spin manifold since $E$ is one.

\begin{lem}
\label{pi1}
Let $N \unlhd F(m_0,b_1, \dots ,b_\nu)$ be normally generated by
$[m_0,m_0^{b_1}], \dots ,$\break $[m_0,m_0^{b_\nu}]$.
Then there is a homomorphism
\[
\varphi_0 \co  \pi_1 \svf \to F(m_0,b_1, \dots ,b_\nu) / N,
\]
such that $\varphi(m) = m_0$ and $\varphi(g_i) = b_i$ for
$i=1,\dots \nu$.
Moreover, this map factors through
$\pi_1 \svf \to \pi_1 W$ induced from the inclusion
$\svf \subset \svfp \hra W$.
\end{lem}

\begin{proof}
To simplify notation we will assume $\nu =1$. The general case can
be worked out analogously.
Let
$B \co  \R^4 \times (0,4) \hra E$
be a parameterization of a neighborhood of $\beta_1$ in
$E$ such that
$B \cap \beta_1 = 0 \times [1,2]$ and
$B \cap K_0 = (0 \times \R^3) \times 1 \cup (\R^3 \times 0) \times 2$.
We may assume that $B \cap K_1 = S \cup T$ where
\[
\begin{array}{lllllll}
   S= &
   &
   (\R &
   \times &
   \R^2 \times 0) &
   \times &
   2
   \\
   T= &
   &
   (0 &
   \times &
   \R^3 - D^3) &
   \times &
   1
   \\
   &
   \cup &
   (0 &
   \times &
   S^2) &
   \times &
   [1,3]
   \\
   &
   \cup &
   (0 &
   \times &
   D^3) &
   \times &
   3 \mbox{ .}
\end{array}
\]
This can in fact be taken as a definition of a finger move.
\begin{figure}[ht]
\vspace{.01in}
\hspace*{\fill}
\relabelbox\small
\epsfysize=2in
\epsfbox{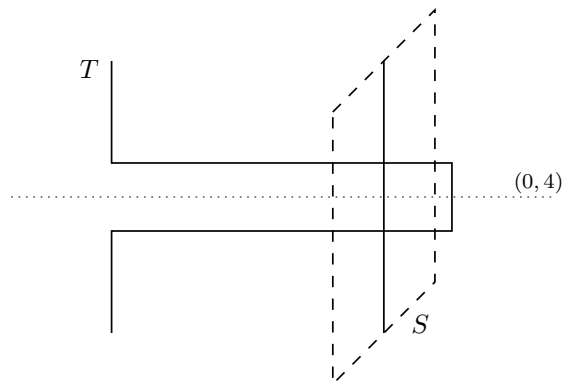}
\relabel {S}{$S$}
\relabel {T}{$T$}
\relabel {0}{\scriptsize$(0,4)$}
\endrelabelbox
\hspace*{\fill}
\vspace{.01in}
\caption{\label{finger} A finger move}
\end{figure}

Let $\bar{B} \cong D^5$ be the compactification of $B$. We have
\[
(E - K_1) = (E-(K_1 \cup B))
          \cup_{\dd \bar{B} - K_1} {(\bar{B}-K_1)}.
\]
Note that
$\pi_1(E - (K_0 \cup B)) \cong
 \pi_1(E - K_0) \cong F(m_0,b_1)$, where $m_0$ is a meridian
to $K_0$ and $b_1$ is an element corresponding to the handle.
$(B-K_1)$ is the complement of two intersecting slice
disks ($S$ and $T$) in the $5$--ball and so by Alexander duality
$H_1 (B-K_1) \cong \Z^2$. (Actually $\pi_1 (B-K_1) \cong \Z^2$
since both slice disks are unknotted.) The corresponding map
$\pi_1 (B-K_1) \to \Z^2$ sends
meridians $m_S$ and $m_T$ of $S$ and $T$
to the generators. $\dd \bar{B} \cap (S \cup T)$ is a trivial
two component link and thus
$\pi_1({\dd \bar{B} - K_1}) \cong F(m_S,m_T)$.
We may assume that
$\pi_1({\dd \bar{B} - K_1}) \to \pi_1 (E-(K_1 \cup B))$ sends the
meridians to $S$ and $T$ to $m_0$ and $m_0^{b_1}$.
Van Kampen's Theorem gives now a map
\[
 \pi_1 (E - K_1) \to F(m_0,b_1)*_{F(m_S,m_T)} \Z^2 \cong
                 F(m_0,b_1) / [m_0,m_0^{b_1}].
\]
The inclusion $\svf \hra (E - K_1)$ sends $m$ to $m_0$
and $g_1$ to $b_1$
(up to some power of $m_0$) and so the inclusions
$\svf \hra W \hra (E - K_1)$ induce the desired map.
\end{proof}

\section{Singular slice disks for embedded links}

Let ${L} \co S^2 \amalg \dots \amalg S^2 \hra S^4$
be an embedded link of $n$
components. Let $X^4$ be the manifold obtained by surgery on $L$.
There is no ambiguity about the framing because $\pi_2 SO(2)=0$.
Note that
$X$ bounds
$D^5 \cup_{{L} \times D^2} \amalg^n D^3 \times D^2$.
The
later manifold possesses a unique spin structure and induces
therefore one on $X$. Let $\pi$ be the fundamental group of $X$
or equally, of the link complement.

Let $m_1, \dots ,m_n$ be a choice of meridians to the
components of
${L}$.
Let $k \geq n+1$, so by Lemma~\ref{nilpotent}
$MF(x_1,\dots,x_n)_k = 0$. The map
$\Phi$ in the following diagram is given by $\Phi(x_i) = m_i$.
We abbreviate $F:=F(x_1,\dots,x_n)$.
\[
\xymatrix
{
 MF  \ar[d]_\cong &
 F  \ar[l] \ar[r]^{\hspace{1mm} \Phi} \ar[d] &
 \pi  \ar[d]
 \\
 MF/MF_{k} &
 F/F_{k}     \ar[r]^\cong \ar[l] &
 \pi / \pi_k
 \\
}
\]
By Alexander Duality $H_2(S^4-{L}) = 0$
and hence $H_2 \pi=0$.
Since
$H_1 \pi = H_1 (S^4 - {L})$ is freely generated by the
meridians,
$\Phi$ induces an isomorphism $H_1 F \cong H_1 \pi$. So by
Stalling's Theorem \cite{S}
the isomorphism in the bottom line of the
above diagram follows. The diagram gives a map
$\Psi \co  \pi \to MF$.
For any $k \geq n+1$ this map factors as
\[
\pi \to F/F_k \to F/F_{n+1} \to MF.
\]
\begin{lem}
\label{boundX}
$X$ bounds a spin manifold $A_0$ with $\pi_1 A_0 \cong MF$
such that
$\Psi$ is the induced map on fundamental groups.
\end{lem}

\begin{proof}
Consider
$\alpha := [X,\sigma,\Psi_0] \in  \Omega_4^{spin}(MF)$.
Here $\sigma$ is the spin
structure on $X$ constructed above and $\Psi_0\co X \to K(MF,1)$
is the map that gives $\Psi$ on the
fundamental group. Since $[X,\sigma]$ is a spin boundary we
have  $\alpha \in {\om}_4^{spin} (MF)$. The
factorization of $\Psi$ implies that $\alpha$ is in the image
of the composition
\[
{\om}_4^{spin} (F/F_k) \to
 {\om}_4^{spin} (F/F_{n+1}) \to
 {\om}_4^{spin} (MF)
\]
for every $k \geq n+1$.
Using $\Omega^{spin}_3(*)=0$,
Theorem~\ref{eventual}
implies now $\alpha = 0$
and hence $X$ bounds a spin
manifold $A_0$
such that
\[
\xymatrix
{
 \pi_1 X \ar[r] \ar[rd]_\Psi &
 \pi_1 A_0 \ar[d]
 \\
 &
 MF
 \\
}
\]
commutes. It is now a standard procedure to do spin structure
preserving surgeries on circles in the interior of
$A_0$ to
obtain $\pi_1 A_0 \cong MF$.
\end{proof}

From  Corollary~\ref{finitely_presented} it follows
that $MF(x_1, \dots ,x_n)$ can be constructed from
$F(x_1, \dots ,x_n)$ by introducing relations
$[x_i,x_i^{h_{i,j}}] = 1$  for $i=1, \dots ,n$ and
$j=1, \dots ,\nu$
and some $h_{i,j} \in F(x_1, \dots ,x_n)$. Fix this $\nu$ and the
$h_{i,j}$ from now on. Let
\[
Y^4 = (S^4 - {L} \times D^2) \cup_{{L} \times S^1}
  \amalg^n \svf.
\]
\begin{figure}[ht]
\vspace{.01in}
\hspace*{\fill}
\epsfysize=3in
\epsfbox{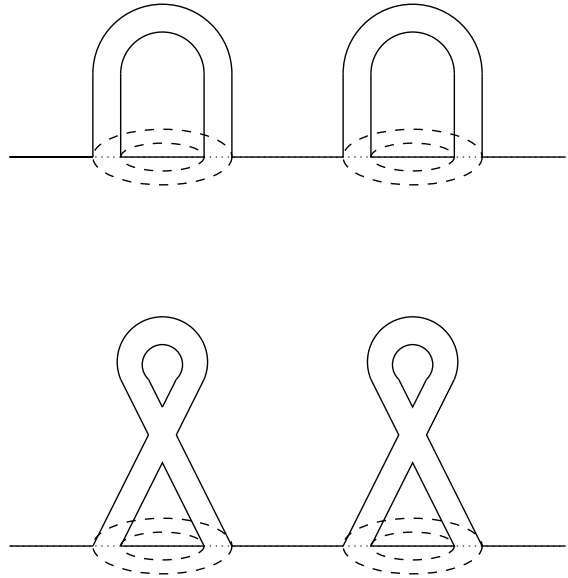}
\hspace*{\fill}
\vspace{.01in}
\caption{\label{xandy} Attaching handles to obtain $X$ and $Y$}
\end{figure}
So $Y$ is the boundary of
$D^5 \cup_{{L} \times D^2} \amalg^n \svh$ (see figure~\ref{xandy}).
We will enumerate the $n$ copies of $\svf$ in $Y$ as
$\svf_1, \dots ,\svf_n$. We will also write
$m,g_1, \dots ,g_n \in \pi_1 \svf$ as
$m_i,g_{i,1}, \dots ,g_{i,\nu}$
when considered as elements of $\pi_1 \svf_i$. By van Kampen's
Theorem we now have
\[
\pi_1 Y = \pi *_{m_1} \pi_1 \svf_1 *_{m_2} \dots
  *_{m_n} \pi_1 \svf_n.
\]
Using $\varphi_0$ from Lemma~\ref{pi1} we can construct  a map
$\varphi \co  \pi_1 Y \to MF$
such that
$\varphi(m_i)=x_i$,
$\varphi(g_{i,j}) = h_{i,j}$ and $\varphi = \Psi$ on $\pi$.

\begin{lem}
\label{boundY}
$Y$ bounds a spin manifold $A$ with $\pi_1 A \cong MF$ such that
$\varphi$ is the induced map on the fundamental group.
\end{lem}

\begin{proof}
Recall that
$X=(S^4 - {L} \times D^2) \cup_{{L} \times S^1}
   \amalg^n D^3 \times S^1$
is  bounded by $A_0$ by Lemma~\ref{boundX}. Let $W_1, \dots ,W_n$
be $n$ copies of $W$, the manifold from Lemma~\ref{pi1}. We may write
$\partial W_i = \svfp_i =
   \svf_i \cup_{S^2 \times S^1} D^3 \times S^1$.
Now we glue $A_0$ to $\amalg_i W_i$ along parts of their boundary to
obtain a manifold
\[
A := A_0 \cup_{\amalg^n D^3 \times S^1} \amalg_i W_i
\]
with boundary $Y$.
%\begin{figure}
%\vspace{6cm}
%\caption{The manifold $A$.}
%\end{figure}
%Using van Kampen's Theorem
\[
\pi_1 A =
  \pi_1 A_0 *_{m_1} \pi_1 W_1 *_{m_2} \dots *_{m_n} \pi_1 W_n.
\]
Using Lemma~\ref{pi1}  and $\pi_1 A_0 \cong MF$
we can find  a map  $\pi_1 A \to MF$ making the triangle
\[
\xymatrix
{
 \pi_1 Y \ar[r] \ar[dr]_\varphi &
 \pi_1 A \ar[d]
 \\
 &
 MF
 \\
}
\]
commute. By the Mayer--Vietoris sequence
\[
H^2 (A;\Z_2) \cong H^2 (A_0;\Z_2) \oplus H^2 (W_1;\Z_2)
  \oplus \dots \oplus H^2 (W_n;\Z_2).
\]
All manifolds on the
right hand side are spin and have therefore vanishing second
Stiefel--Whitney class. Hence  $A$  is also spin. Again we
can do spin structure preserving surgeries in the interior of
$A$ to obtain $\pi_1 A \cong MF$.
\end{proof}

\begin{thm}
$L$ bounds disjointly immersed slice disks in $D^5$.
\end{thm}

\begin{proof}
Let $A$ be the spin manifold obtained in Lemma~\ref{boundY}.
Recall that
\[
\partial A = Y =
 (S^4 - {L} \times D^2) \cup_{{L} \times S^1}
 \amalg_i \svf_i
\]
Now let $\Sigma H := \svh_1 \amalg \dots \amalg \svh_n$ be the
disjoint union
of $n$ copies of $\svh$. Then $\partial \Sigma H$ and $\partial A$
both contain the disjoint union of $n$ copies of $\svf$ which we will
denote by $\Sigma F$. Now glue $\Sigma H$ to $A$ along $\Sigma F$ to
obtain
\[
D := A \cup_{\Sigma F} \Sigma H.
\]
Then
$\partial D =
 (S^4 - {L} \times D^2) \cup_{{L} \times S^1}
 ({L} \times D^2) = S^4$
and ${L}$ bounds disjointly immersed slice disks in
$\Sigma H \subset D$. We will modify $D$ to the $5$--ball
without destroying the immersed slice disks. In order to do so
it is enough to make $D$ 2--connected (using
Poincar{\'e} duality and the h--cobordism Theorem). Again by
van Kampen's Theorem
\[
\pi_1 D = \pi_1 A *_{\pi_1 \svf_1} \pi_1 \svh_1
  *_{\pi_1 \svf_2} \dots *_{\pi_1 \svf_n} \pi_1 \svh_n
\]
and this group is normally generated by the
meridians $m_i$
(the $g_{i,j}$ are identified with products of meridians). But
the meridians bound $2$--disks in $\Sigma H$ and so
$\pi_1 D=1$. To kill  $\pi_2 D$ we will  do surgery on
$2$--spheres in $A$ (and hence away from the immersed slice disks).

The horizontal line of the  diagram below is part of
a Mayer--Vietoris sequence. By Lemma~\ref{pi1}
$H_1 \Sigma F \cong H_1 A \oplus H_1 \Sigma H$, and this
explains that $H_2 D \to H_1 \Sigma F$ below is trivial.
\[
\xymatrix
{
 &
 \pi_2 A \ar[r] \ar[d] &
 \pi_2 D \ar[d]_\cong
 \\
 H_2 \Sigma F \ar[r] \ar[rd] &
 H_2 A \oplus H_2 \Sigma H \ar[d] \ar[r] &
 H_2 D \ar[r]^0 &
 H_1 \Sigma F  \ar[r]^{\cong} &
H_1 A \oplus H_1 \Sigma H
 \\
 &
 H_2 K(MF,1)
}
\]
Now according to Lemma~\ref{onto}
$H_2 \Sigma F$  maps onto $H_2 \Sigma H$
and the Clifford tori represent elements in the kernel of this map.
But the Clifford tori also
generate $H_2 K(MF,1)$ \cite{FT}, and
so $H_2 \Sigma F$ maps also onto $H_2 K(MF,1)$.
The exact sequence
\[
\pi_2 A \to H_2 A \to H_2 K(MF,1) \to 0
\]
proves now that $\pi_2 A$ maps onto
$\pi_2 D \cong H_2 D$.

It is a classical result of Milnor and Kervaire
in \cite{KM}
that $D$ can be changed to the $5$--ball by a sequence of surgeries
on classes in $\pi_2 D$.  (Their original
formulation considers framed manifolds, but here a spin structure
is enough to guarantee trivial normal bundle on  $2$--spheres.)
We just saw that $\pi_2 A$ maps onto $\pi_2 D$, so we can represent
these classes by  $2$--spheres in $A$
(and by general position embedding  these
spheres comes for free). So we can do surgery to $A$ and change
$D$ to the $5$--ball.
\end{proof}

\section{Link concordance implies link homotopy} \label{sec:habil}

Let $M^3$ be a 3--manifold with boundary
$\dd M = \dd_0 M \amalg \dd_1 M$. Consider a generic immersion
$f \co  M \imra N^4 \x I$, where $N$ is a closed $4$--manifold and $f^{-1}(N^4
\x\{i\}) =\dd_i M$  for $i=0,1$. The genericity of $f$ implies that the
double point set $S(f)$ is a 1--dimensional manifold with boundary in $\dd
M$. We say that
$f$ is a {\em Morse immersion} if
$p_2\circ f$ is a generic Morse function on $M$ and so is the
restriction to $S(f)$. Recall that Morse functions are generic iff all
critical values are distinct. For a Morse immersion we also assume that the
critical values of
$p_2\circ f$ and $p_2\circ f|_{S(f)}$ are distinct. Regular homotopies
of a closed surface $F^2$ give examples of Morse immersions
$F^2\x I\to N^4\x I$ without critical points for $p_2\circ f$.

Let $g \co  F^2 \imra N^4$ be
a generic immersion. We describe six ways to
construct a Morse immersion
$f \co  M \imra N^4 \x I$ with $\dd_0 M = F$ and $f_0 = g$
using a {\em guiding map} $\alpha$.
To simplify the discussion we will ignore framing data on
$\alpha$ in these examples. They turn out to be irrelevant for our discussion.
\begin{itemize}
\item[($h1$)] Let $\alpha \co  I \hra N^4$ be guiding arc for a
            fingermove on $g$. The fingermove along
            $\alpha$ gives a
            Morse immersion  with $M = F \x I$ and
            $f_1 = $ result of a fingermove on $g$.
\item[($h2$)] Let $\alpha \co  D^2 \hra N^4$ be a Whitney disk
            for $g$.
            The Whitney move along $\alpha$ gives a
            Morse immersion with $M = F \x I$ and
            $f_1 = $ result of a Whitney move on $g$.
\item[($b0$)] Let $\alpha \co  D^0  \hra N^4$ be a point
            in the complement of $g$. Let
            $h \co  D^3 \hra N^4$ be an embedding into a regular
            neighborhood of $\alpha$. This data gives rise to
            Morse immersion with $M = (F \x I) \amalg D^3$,
            $\dd_1 M = F \amalg S^2$ and
            $f_1 = g \amalg h|_{S^2}$.
\item[($b1$)] Let $\alpha \co  D^1 \hra N^4$ be an embedded arc such that
            $\alpha(D^1) \cap g(F) = \alpha(\dd D^1)$ and that $\alpha$
            misses all double points of $g$. We can use $\alpha$ to
            do ambient surgery on $g$ to obtain a Morse immersion with
            $M = F \x I \cup$ 1--handle and
            $\dd_1 M =$ result of a surgery on an $S^0 \subset F$.
\item[($b2$)] Let $\alpha \co  D^2 \hra N^4$ be an embedding  such that
            $\alpha(D^2) \cap g(F) = \alpha(\dd D^2)$
            and that $\alpha$
            misses all double points of $g$. We can use $\alpha$ to
            do ambient surgery on $g$ to obtain a Morse immersion with
            $M = F \x I \cup$ 2--handle and
            $\dd_1 M=$ result of surgery on an $S^1 \subset F$.
\item[($b3$)] Let $\alpha \co  D^3 \hra N^4$ be an embedding
            such that $\alpha|_{S^2}$ is the restriction of
            $g$ to one component of $F$ and the interior of $\alpha(D^3)$
            misses $g$. This gives a Morse immersion with
            $M = F \x I \amalg D^3$, $\dd_1 M \amalg S^2 = F$
            and $f_1= g|_{F-S^2}$.
\end{itemize}
 We will call such Morse immersions  {\em elementary} and
will refer to $h1,\dots,b3$
as the {\em type} of an elementary Morse immersion. For type
$h1$ and $h2$ the Morse function $p_2\circ f$
has no critical points whereas the restriction to $S(f)$ has a
minimum for $h1$ and a maximum for $h2$.
The corresponding regular homotopies are the finger move (for $h1$) and the
Whitney move (for $h2$).

In the other  cases the
restriction to $S(f)$ has no critical points but $p_2\circ f$ has exactly one
critical point of index $i$ for $bi$. These describe the ambient handle
decomposition of the $3$--manifold $M$.

 By turning the interval $I$ upside
down, every elementary Morse immersion
$f \co  M \imra N^4 \x I$ can also be reconstructed from
$f_1 \co  \dd_1 M \imra N^4 \x \{1\}$ and
an attaching map $\beta$ into $N^4 \x \{1\}$.
This changes the types
as follows :  $h1 \lra h2, b0 \lra b3, b1 \lra b2$.

Note that in
each  elementary Morse immersion everything happens in a regular neighborhood
of $\alpha \x I \subset N^4 \x I$. In particular,
the complement of this neighborhood
is a product (with respect to the
product structure on $N^4 \x I$).

The next lemma
states that Morse immersions are generic and that every Morse immersion
can be constructed from elementary ones. We believe this general position
result is  well known, details can be found in
\cite{habil} (for all dimensions). In the following, an isotopy of
immersions is the conjugation of the map by diffeotopies of range and
domain, just as in the stability condition for smooth mappings \cite[Section
III]{GG}.

\begin{lem}[General Position]
\label{zerlegung}
Let $f \co  M^3 \imra N^4 \x I$ be a generic immersion.
It is isotopic (rel $\dd$) to a Morse immersion $g \co  M \imra N^4
\x I$. After a further isotopy, we may assume that there are values
$0=a_0 < a_1 < \dots < a_r=1$ such that the restrictions $g^i$ of $g$
to $(p_2\circ g)^{-1}[a_i,a_{i+1}]$ are elementary Morse immersions (after
      rescaling $[a_i,a_{i+1}]$ to $[0,1]$).
\end{lem}
Note that the index $i$ of the elementary Morse immersion $bi$ is the
dimension of the guiding map (or descending manifold) $ \alpha$ and thus
it is consistent to give $hi$ index $i$ as well. This orders all the
critical points of a Morse immersion.

\begin{lem}[Reordering]
\label{tausch}
In the notation of Lemma~\ref{zerlegung}, if index $g^{i+1} \leq$ index
$g^i$, then g is isotopic to a  Morse immersion $h\co  M^3 \imra N^4 \x I$ such
that each $h^i$ is elementary and
\begin{itemize}
\item[\rm(i)]
     $g^j = h^j$ for $j \neq i,i+1$;
\item[\rm(ii)]
     type $h^{i}=$ type $g^{i+1}$ and
     type $h^{i+1}=$ type $g^{i}$.
\end{itemize}
\end{lem}

\begin{proof}
Let $\alpha_{i+1}$ (resp.\ $\beta_{i}$) be the attaching maps
into $N^4 \x \{a_{i+1}\}$
that guide the construction of $g^{i+1}$ (resp. $g^{i}$) upwards (resp.\
downwards) from
$g^{i+1}_0 = g^{i}_1$. Now index $g^{i+1} \leq$ index $g^i$
implies
$$
\dim \alpha_{i+1} + \dim \beta_i = \ind g^{i+1} + (3-\ind g^i)\leq 3.
$$
By general position we may assume
that $\alpha_{i+1}$ and $\beta_i$ are disjoint in the level $N^4 \x
\{a_{i+1}\}$.  But this means that $g^i$ and
$g^{i+1}$ are independent in the following sense:
$\alpha_{i+1}$ can be
pushed down to an embedding $\alpha'_{i}$
into $N^4 \x \{a_i\}$ that does
not intersect the attaching map $\alpha_{i}$ for
$f^i$ into $N^4 \x \{a_i\}$. Now first $h^{i}$ is constructed from
$\alpha'_{i}$ and then $h^{i+1}$ from $\alpha_{i}$.
\end{proof}

Recall that a link map is a continuous map which keeps distinct
components disjoint. A link concordance is a link map $f\co F \x I \to N\x I$
such that $f^{-1}(N\x \{i\})=F \x \{i\}$ for $i=0,1$. Finally,
 a link homotopy is a homotopy through link maps.
We make a very simple, but useful
observation: If a link concordance $f$ is a product on all but one
component, then $p_1 \circ f$ is a link homotopy. By applying the above
general position and reordering lemmas, we will repeatedly be able to
apply this observation to prove the main result of this section:

\begin{thm}[Link concordance implies link homotopy] \label{four}
If $f\co F^2 \x I \to N^4\x I$ is a link concordance then there is a link
homotopy $h\co F^2 \x I\to N^4$ such that $f_i=h_i$ for $i=0,1$.
\end{thm}

\begin{proof}
Using Thom's jet transversality theorem \cite[Sections II.4--5]{GG} we can assume
that $f$ is a generic immersion except for a finite number of cross caps
\cite[page 179]{GG}. As in Whitney's original immersion argument, these
cross caps can be pushed off the, say, lower boundary $N\x \{0\}$. This
changes the lower boundary $f_0$ by cusp homotopies which we may assume
are small enough such that
they don't change the link homotopy class of $f_0$.

Let $g\co  F^2\x I\imra N^4\x I$  be a Morse immersion
satisfying the conclusion from Lemma~\ref{zerlegung}.
If all  the $g^i$ are of types
$h1$ and $h2$ then $p_1\circ g$ is in fact a link homotopy and we
are done. Let $C_1 ,\dots, C_n$ be the components of $F^2 \x I$.
We know that $C_1 \cong M^2 \x I$. Consequently,
all critical points of index $0$ of $p_2 \circ g$ on $C_1$ can be
canceled (by Morse cancellation) by critical points of index $1$.
By applying Lemma~\ref{tausch} several times to move other critical points up
respectively down, we may assume that the corresponding elementary  Morse
immersions for the canceling 0-- and 1--handles of $C_1$ are consecutive. Let's
say they are
$g^{k},g^{k+1},\dots,g^{l}$ and let $W :=
g^{-1}([a_k,a_{l+1}])$. Then  $W \cap C_1 \cong M^2 \x I$ and
$g|_W \co  W \imra N^4 \x [a_k,a_{l+1}]$ is
a product on $W - C_1$.
Hence the observation from above applies and we can replace
$g^{k},g^{k+1},\dots,g^{l}$ by elementary Morse immersions of
type $h1$ and $h2$.
More precisely, we replace
$$
g(m,t)=(g_1(m,t),g_2(m,t)), \quad t\in [a_k,a_{l+1}]
$$
by $(g_1(m,t),t)$ which doesn't change anything away from $W$ and removes the
critical points on $W$. Finally, we can make this map generic, producing a
Morse immersion with singularities of types $h1$ and $h2$ only.

We can apply the same procedure to all components
and all $g^i$ of type $b0$ and (by symmetry) $b3$.
Hence we may now assume that all the $g^i$
are of type $h1\ $,$h2\ $,$b1$
and $b2$. Again by Lemma~\ref{tausch} we can order the types of
the $g^i$, such that $g^{k},g^{k+1},\dots,g^{l}$
are consecutive
elementary Morse immersions that induce critical points on $C_1$. We still
know $C_1 \cong  M^2 \x I$ and the same argument as before
allows us
to replace $g^{k},g^{k+1},\dots g^{l}$
by elementary Morse immersions of type $h1$ and $h2$. Repeating this
procedure on $C_2,\dots,C_n$ finishes the proof.
\end{proof}

\begin{rem} \label{one}
The first part of the argument which makes $f$ a generic
immersion is not really necessary to prove our Theorem~\ref{main} because
the output of the previous sections is already a generic immersion by
construction. We just included this step for completeness.
\end{rem}
\begin{rem}
The above argument does not explicitly mention handle cancellation except in
the elementary case involving $0$--handles, ie, critical points of type $b0$
(and $b3$). This is done on purpose because only in this simplest case can one
avoid to use {\em ambient handle slides} in order to obtain canceling pairs of
handles. However, if one is willing to introduce a gradient-like vector field
more explicitly into the discussion, then handle slides are well defined
and can
be in fact done ambiently. This is used by Rourke in \cite{R} to do all the
steps of the proof of the h--cobordism theorem ambiently, obtaining a proof of
Hudson's  ``concordance implies isotopy'' in codimension $\geq 3$ (this
restriction is explained in the next remark). Rourke's argument works for
embeddings and is generalized to generic immersions in \cite{habil},
 where the multiple point stratification on the range of a generic immersion
is considered and  stratified versions of Morse
functions, gradient-like vector fields and handles are introduced.
\end{rem}
\begin{rem}
The reason why Hudson's theorem does not work in codimension~$2$ can already
be seen for knotted arcs in $D^2 \times I$. More precisely, if one tries to
ambiently cancel two handles which do cancel abstractly, one has to push the
core $C$ of the higher handle together with the cocore $Q$ of the lower handle
into the middle level (this uses the gradient-like vector field). A dimension
count shows that they may be assumed disjoint in codimension $\geq 3$ and
hence can be canceled ambiently. However, in codimension~$2$ one can
only assume that
$C$ and $Q$ intersect transversaly in a finite number of points.
This produces basically all possible knotting phenomena.
But since $C$
and $Q$ are part of the same component, it is not a problem in the link
homotopy world: One just maps the abstract cancellation forward into a
neighborhood of $C\cup Q$, producing a link homotopy (rel.\ boundary) to the
situation where $C$ and $Q$ have canceled. By thickening into the dual
dimensions of $C$ and $Q$ one sees that each intersection point in
$C\cap Q$ contributes to a small sphere of self-intersections on the relevant
component. In the dimension range of Theorem~\ref{four} $C$ and $Q$ are
$2$--disks and the new self-intersections are circles of double points. Hence in
addition to removing a single pair of critical points of types $b1$ and
$b2$ the
procedure introduces one local pair of critical points of type
$h1$ and $h2$ for each point in $C\cap Q$. This together with
 Remark~\ref{one} explains how the proof of Theorem~\ref{main} could
be given completely in the category of generic immersions.
\end{rem}
\begin{rem}
In the argument of \cite{R} explained above one has to assume that the
dimension
range is such that the Whitney trick can be applied. Since this excludes some
very  interesting low dimensions, the point of view in
\cite{habil} is different. There ambient Cerf theory is used to give a proof
of Hudson's theorem in all dimensions. This method also avoids the explicit
mentioning of handle slides but nevertheless the stratified versions of
gradient
like vector fields are essential.
\end{rem}
\begin{rem}
In \cite{habil} it is proven that a link concordance $M^m
\x I\imra N^n \x I$, which is an immersion, is homotopic (rel. boundary) to
a link
homotopy if
$m\leq n-2$. The higher dimensional analogues of cross caps, ie, points
where the map is not an immersion, are not discussed. However, by
Hirsch--Smale immersion theory, one still obtains that link concordance
implies link homotopy for maps $S^{m_1}\amalg \dots\amalg S^{m_r}\to S^n$
if $m_i\leq n-2$.
\end{rem}

\end{document}